%% file: qithree.tex
\documentclass[12pt,draft]{article}

\usepackage{amsmath,amsthm,latexsym}

\def\var{ma\-ni\-fold}
\def\svar{sub\-ma\-ni\-fold}
\def\Int{\operatorname{Int}}
\def\diam{\operatorname{diam}}
\def\inj{\operatorname{inj}}
\def\capac{\operatorname{cap}}
\def\bord{\partial}
\def\cqfd{\ensuremath{\Box}}

\def\gen{ge\-ne\-ra\-ted}

\def\qi{quasi-iso\-me\-tric}

\def\qy{quasi-iso\-me\-try}
\def\geod{geo\-de\-sic}
\def\qg{quasi-geo\-de\-sic}
\def\bg{boun\-ded geo\-me\-try}

\def\homeoic{ho\-meo\-mor\-phic}
\def\diffeo{dif\-feo\-mor\-phism}

\def\riem{Rie\-man\-nian}
\def\virt{vir\-tu\-ally}
\def\isop{i\-so\-pe\-ri\-me\-tric}
\def\hyp{\mathrm{hyp}}
\def\vol{\mathrm{vol}}
\def\sing{\operatorname{sing}}

\let\bydef\emph

\newtheorem{theo}{Theorem}[section]
\newtheorem{corol}[theo]{Corollary}
\newtheorem{prop}[theo]{Proposition}
\newtheorem{lem}[theo]{Lemma}
\newtheorem{slem}[theo]{Sublemma}

\theoremstyle{definition}
\newtheorem*{defi}{Definition}
\newtheorem*{claim}{Claim}

\newtheorem*{rem}{Remark}
\newtheorem*{Notation}{Notation}
\newtheorem*{ackno}{Acknowledgments}

\def\Rr{\mathbf{R}}

\def\Zz{\mathbf{Z}}
\def\Nn{\mathbf{N}}

\def\H{\mathbf{H}}
\def\E{\mathbf{E}}
\def\Ss{\mathbf{S}}
\def\slr{\widetilde{\mathrm{SL}_2 (\Rr)}}

\def\color[#1]#2{}

\title{Large-scale conformal rigidity in dimension three}
\author{Sylvain Maillot\thanks{Partially supported by a CRM-CIRGET
postdoctoral fellowship.}\\ \texttt{maillot@math.uqam.ca}}
\date{June 19, 2002}

\begin{document}

\maketitle

\section{Introduction}
Let $X$ be a noncompact smooth \var. Two complete \riem\ metrics on
$X$ in a given conformal class can have very different asymptotic
geometries. For instance, starting with the Euclidean plane $\E^2$
with polar coordinates $(r,\theta)$, multiplying the Euclidean metric
$g_0=dr^2 + r^2 d\theta^2$ by a function equal to $1/r^2$ outside a
compact neighborhood of the origin, one obtains a complete
\riem\ metric $g$ which is \qi\ to a half-line. The metrics $g$ and
$g_0$ are in the same conformal class, but they are not \qi. In fact,
they have different asymptotic dimensions.

One can ask whether something similar can happen to two finitely generated
groups $\Gamma$ and $\Gamma_0$: can they be `coarsely quasi-conformal'
in some sense and yet not \qi?

To give a precise meaning to this question, we will choose Riemannian \var s
$X_0,X$ that are geometric models for $\Gamma_0,\Gamma$ and consider
conformal mappings between $X_0$ and $X$.
First let us recall some definitions. Two metric spaces $(X_1,d_1)$
and $(X_2,d_2)$ are \bydef{\qi} if there are constants $\lambda\ge1$
and $C\ge0$ and a map $f:X_1\to X_2$ satisfying~:
\begin{gather*}
\lambda^{-1}\, d_1(x,x')-C \le d_2(f(x),f(x'))\le
\lambda\, d_1(x,x')+C\qquad \forall x,x'\in X_1\\
\forall y\in X_2,\quad \exists x\in X_1, \quad d(f(x),y) \le C.
\end{gather*}

If a group $\Gamma$ is finitely \gen\ by a subset $S$, one can make
$\Gamma$ into a metric space by means of the word metric associated to
$S$. The quasi-isometry class of this metric does not depend on the
choice of $S$. Thus one can omit the mention of $S$ when discussing
quasi-isometries between groups and metric spaces.

Let $\Gamma$ be a group and $X$ a metric space.
A \bydef{geometric action} of $\Gamma$ on $X$ is a proper, cocompact
action of $\Gamma$ on $X$ by isometries.
We are interested in the case where $X$ is a Riemannian
\var. The Riemannian \var s we consider will always be complete and of
\bg. For the purpose of this paper, it is convenient to take the
following definition: a \riem\ \var\ $X$ has \bydef{\bg} if there is a 
number $\epsilon>0$ and a compact \riem\ \var\ $Y$ such that all balls
of radius $\epsilon$ in $X$ can be isometrically embedded into $Y$.
It is \bydef{conformally flat} if every point has a neighborhood
conformal to a ball in Euclidean space.

A special case of the ``fundamental observation of geometric group theory''
(\cite[Proposition 8.19]{bh:cat}) is that
if a group $\Gamma$ acts geometrically on a complete \riem\ \var\ $X$,
then $\Gamma$ is finitely \gen\ and \qi\ to $X$. Furthermore, such a
$X$ is clearly of \bg.

We can now state our problem more formally: let $\Gamma_0$ be a
finitely \gen\ group acting geometrically on a complete \riem\ \var\ $X_0$. 
The `coarse quasi-conformal' deformations we are looking for
are pairs $(\Gamma,X)$ where $\Gamma$ is a finitely \gen\ group and $X$ a
complete \riem\ \var\ of \bg\ \qi\ to $\Gamma$ and conformal to $X_0$.
We regard a deformation as trivial if all spaces involved are
\qi. Since the existence of nontrivial deformations depends only on
$X_0$ and not on $\Gamma_0$, we make the following definition.

\begin{defi}
  Let $X_0$ be a complete Riemannian \var\ which admits a geometric
  group action. We say that $X_0$ is \bydef{large-scale conformally
  rigid} if every finitely \gen\ group \qi\ to a complete Riemannian
  \var\ of \bg\ conformal to $X_0$ is in fact \qi\ to $X_0$.
\end{defi}

Our main result is the following rigidity theorem:
\begin{theo}\label{prince}
  Let $X_0$ be a complete Riemannian $3$-\var\ which admits a geometric
  group action. Assume that $X_0$ is conformally flat and
  \homeoic\ to $\Rr^3$. Then $X_0$ is large-scale conformally rigid.
\end{theo}

(In fact we prove a little more; see the last section for a discussion
of this.)

Theorem~\ref{prince} applies to three among Thurston's eight
$3$-dimensional geometries, namely $\E^3$, $\H^3$ and $\H^2\times\Rr$.
Groups which are \qi\ to those geometries are
known~\cite{gromov:growth, cc:hyper, rieffel:groups}. In particular,
we obtain the following characterization of groups acting
geometrically on $\E^3$ and $\H^3$:

\begin{corol}
Let $\Gamma$ be a finitely \gen\ group. Then $\Gamma$ admits a
geometric action on $\E^3$ (resp.~$\H^3$)
if and only if it is \qi\ to some complete Riemannian \var\ of \bg\ conformal
to $\E^3$ (resp.~$\H^3$).
\end{corol}

The proof of Theorem~\ref{prince} splits into two cases, according to
whether the conformal structure of $X_0$ is `parabolic' or
`hyperbolic'. The necessary background is reviewed in Section~\ref{sec
parabole}. The parabolic case of Theorem~\ref{prince} is proved in
Section~\ref{sec cas parabole} using results from Section~\ref{sec
parabole} and arguments from coarse topology. The hyperbolic case is
tackled in Section~\ref{sec cas hyperbole}. Various remarks on
generalizations and open questions are gathered in Section~\ref{sec final}.

\begin{ackno}
The author wishes to thank C.~Pittet and T.~Delzant for conversations
related to this work.
\end{ackno}

\begin{Notation}
When $A$ is a subset of a finitely \gen\ group or a \riem\ \var, we denote by
$|A|$ its ``volume'', i.e.~: if $A$ is finite (resp.~a curve,
resp.~a surface, resp.~a domain with nonempty interior), $|A|$
denotes the cardinal of $A$ (resp.~its length, resp.~its area, resp.~its
volume.)

We systematically denote by $d$ the distance function of a metric space.
A metric ball (resp.~sphere) arounf a point $x$ of radius $r$ is
denoted by $B(x,r)$ (resp.~$S(x,r)$).
\end{Notation}

\section{Discrete groups and $p$-parabolicity}\label{sec parabole}
\subsection{A review of $p$-parabolicity}
Throughout this subsection we fix an integer $p\ge 2$.

\begin{defi}
Let $X$ be a \riem\ \var. The $p$-\bydef{capacity} of a compact subset
$K\subset X$ is defined by $$\capac_p(K) = \inf_u \int_X |\nabla
u|^p\, d\vol,$$ where the infimum is taken over all compactly
supported smooth functions $u$ such that $u(x)\ge 1$ for every $x\in K$.

The \var\ $X$ is $p$-\bydef{parabolic} if $\capac_p(K)=0$ for every compact
$K\subset X$. Otherwise it is $p$-\bydef{hyperbolic}.
\end{defi}

Our interest in this notion comes from the following fact:
if $X_0$ and $X$ are \riem\ $p$-\var s conformal to
each other, then $X_0$ is $p$-parabolic iff $X$ is $p$-parabolic.
Moreover, $p$-pa\-ra\-bo\-li\-ci\-ty is a \qy\ invariant of complete
\riem\ \var s of \bg~\cite{kanai:parabolicity,holo:rough}. However we
will not need this result: we are interested in \var s that are \qi\
to groups, and in this case the \qy\ invariance is a consequence of a
characterization in terms of growth functions and \isop\ inequalities
(Theorem~\ref{varogromov} below.)

\begin{rem}
$2$-parabolicity is equivalent to the recurrence of the Brownian motion, or
to the existence of a Green function for the Laplace-Beltrami
operator. The relevance of these ideas to large-scale conformal
rigidity in dimension~2 was observed by G.~Mess~\cite{mess:seifert}.
For more information and references for the general case,
see~\cite{grigoryan:capacities}.
\end{rem}

\subsection{Growth and isoperimetry}

Let $\Gamma$ be a finitely \gen\ group and $S$ a finite generating
subset of $\Gamma$. For all $\Omega\subset\Gamma$ we set
$$\bord\Omega:=\{g\in\Omega \mid \exists
g'\in\Gamma-\Omega,\ d_S(g,g')=1\}.$$

In the following definitions, $\Gamma$ is a finitely \gen\ group and $X$ a
\riem\ \var.

The \bydef{growth function} of $\Gamma$ (resp.~$X$) is the function
$r\mapsto |B(*,r)|$, where $*$ is a basepoint. We say that $\Gamma$
(resp.~$X$) has \bydef{superpolynomial growth} if for each $D>0$ there
exists $C_D>0$ such that $|B(*,r)|\ge C_D n^D$ for all $r$. We say
that $\Gamma$ (resp.~$X$) has \bydef{polynomial growth of exponent}
$D\in\Nn$ if there is $C>0$ such that $C^{-1} r^D\le |B(*,r)| \le C
r^D$ for all $r$.

An \bydef{\isop\ function} for $\Gamma$ (resp.~$X$) 
is a function $I:[0,+\infty) \to [0,+\infty)$ such that the inequality
$$I(|\Omega|) \le |\bord\Omega|$$
holds for every finite subset of $\Gamma$ (resp.~every bounded domain in
$X$ with sufficiently smooth boundary).

The \bydef{\isop\ dimension} of $\Gamma$ (resp.~$X$)
is the supremum of the set of
$D\ge 0$ for which there is a constant $C>0$ such that the
function $v\mapsto C v^{(D-1)/D}$ is an \isop\ function.

A theorem of Gromov~\cite{gromov:growth} says that the growth function
of a group is either superpolynomial or polynomial; in the latter
case, the group is virtually nilpotent and the exponent of growth can
be computed from the ranks of quotients in the lower central
series~\cite{bass:degree}.

The \isop\ dimension and the asymptotic behavior of the growth function of
$\Gamma$ do not depend on the choice of the generating set $S$; in
fact they are \qy\ invariants of groups and complete \riem\ \var s of
\bg~\cite{kanai:rough}.

The following theorem follows from various results scattered in the
literature and does not seem to have been stated before in this
generality. Since we think it is of independent interest, we give a
more complete statement than we shall actually need.
The main ingredients are due to Gromov and Varopoulos.

\begin{theo}\label{varogromov}
Let $\Gamma$ be a finitely \gen\ group and $X$ a complete noncompact
\riem\ \var\ of \bg\ \qi\ to $\Gamma$. The following are equivalent:
\begin{enumerate}
\item $X$ is $p$-parabolic;
\item The \isop\ dimension of $\Gamma$ (or $X$) is at most $p$;
\item $\Gamma$ is \virt\ nilpotent of growth exponent at most $p$.
\end{enumerate}
\end{theo}

\begin{proof}
We begin by proving that (i) implies (ii). By~\cite[section
3]{grigoryan:capacities}, the $p$-parabolicity of $X$ implies that for
every isoperimetric function $I$ the following holds:
$$\int^\infty \, \frac{dv}{I^{\frac{p}{p-1}}(v)} = \infty.$$
Assuming by contradiction that there exists $D>p$ such that
$v\mapsto C v^{(D-1)/D}$ is an \isop\ function for $X$, we get:
$$\int^\infty \,
\frac{dv}{C^{\frac{p}{p-1}}v^{\frac{(D-1)p}{D(p-1)}}}= \infty.$$
This contradicts the fact that $D>p$.

Let us turn to the proof that (ii) implies (iii). If (iii) does not
hold, then Gromov's Theorem~\cite{gromov:growth} implies that the
growth of $\Gamma$ is superpolynomial or polynomial of exponent at
least $p+1$. By Varopoulos's inequality~\cite[Thm 1]{cs:isop}, the
\isop\ dimension of $\Gamma$ is at least $p+1$, so (ii) does not hold.

Finally we show (iii) implies (i). If $\Gamma$ has polynomial growth
of exponent at most $p$, then the same is true for
$X$, i.e.~for all $x\in X$ there is a constant $C$ such that
$|B(x,r)| \le C r^p$ for large $r$.
By~\cite[section 3]{grigoryan:capacities}, to prove $p$-parabolicity
it is enough to check that $$\int^\infty \, \left(\frac
r{|B(x,r)|}\right) ^{1/(p-1)} \,dt = \infty.$$
From the upper bound on $|B(x,r)|$ we deduce $$\frac r{|B(x,r)|} \ge
\frac 1{C^p r^{p-1}},$$ which implies the divergence of the above integral.
\end{proof}

We shall be mostly interested in the case $p=3$, so we state for future
reference a corollary to Theorem~\ref{varogromov}:
\begin{corol}\label{parabolic three}
Let $X$ be a noncompact complete \riem\ \var\ of \bg\ and $\Gamma$ a
finitely \gen\ group \qi\ to $X$. Then $X$ is $3$-parabolic if and
only if $\Gamma$ is \virt\ $\Zz^n$ with $1\le n\le 3$.
\end{corol}

\begin{proof}
The ``if'' part follows from the implication (iii) $\implies$ (i) in
Theorem~\ref{varogromov}. The ``only if'' part follows from the
implication (i) $\implies$ (iii) in the same theorem plus the formula
for the exponent of growth of a nilpotent group~\cite{bass:degree}.
\end{proof}

\section{The 3-parabolic case}\label{sec cas parabole}
The goal of this section is to prove Theorem~\ref{prince} in the case
where $X_0$ is $3$-parabolic. In fact we will prove a stronger
statement:
\begin{theo}\label{parabolic case}
Let $X,X_0$ be complete \riem\ \var s of \bg\ \homeoic\ to $\Rr^3$ and
conformal to each other. Let $\Gamma_0,\Gamma$ be finitely \gen\ groups
\qi\ to resp.~$X,X_0$. If $X_0$ is $3$-parabolic, then both
$\Gamma_0$ and $\Gamma$ are \virt\ $\Zz^3$. In particular, they are \qi\ to
each other.
\end{theo}

Our main tools are Corollary~\ref{parabolic three} and the following
topological rigidity result:
\begin{prop}\label{topo rigidity}
Let $Y$ be a complete \riem\ \var\ of \bg\ \homeoic\ to $\Rr^3$. Then
$Y$ is not \qi\ to $\E^2$.
\end{prop}

\begin{proof}[Proof of Theorem~\ref{parabolic case} assuming
Proposition~\ref{topo rigidity}]
By Corollary~\ref{parabolic three}, $\Gamma_0$ is virtually $\Zz$,
$\Zz^2$ or $\Zz^3$. We must rule out the first two cases. If
$\Gamma_0$ were virtually $\Zz$, it would have two ends. This
contradicts the hypothesis that $X_0$ is \homeoic\ to $\Rr^3$, because
the number of ends is a \qy\ invariant for groups and complete
Riemannian \var s (cf.~\cite[Proposition 8.29]{bh:cat}). The
possibility that $\Gamma_0$ be virtually $\Zz^2$ is prohibited by
Theorem~\ref{topo rigidity}. Hence $\Gamma_0$ is virtually $\Zz^3$.

Since $3$-parabolicity is conformally invariant in dimension $3$, the
same arguments apply to $X$, so $\Gamma$ is virtually $\Zz^3$. In
paticular, $\Gamma,\Gamma_0,X,X_0$ are all \qi\ to $\E^3$.
\end{proof}

The end of this section is devoted to the proof of
Proposition~\ref{topo rigidity}. We first give the idea of the proof,
which is fairly simple.

Seeking a contradiction, we assume that there is a \qy\ $f:\E^2
\to Y$ and fix a coarse inverse $\bar f:Y\to \E^2$. We want to exploit
the fact that $Y$ is simply-connected at infinity and $\E^2$ is
not. With this in mind, take a large round circle $c_1$ in $\E^2$. Its image
by $f$ is a quasi-circle in $Y$; since $Y$ is a \geod\ space, we can
approximate it by a true topological circle $c_2$. Now $Y$ is
simply-connected at infinity, so $c_2$ can be filled by a disc $D_2$
near infinity. Since $\E^2$ is uniformly simply-connected, the
quasi-disc $\bar f(D_2)$ can be approximated by a topological disc,
which stays near infinity and is homotopic near infinity to
$c_1$. This contradicts the fact that $\E^2$ is not simply-connected
at infinity.

Note that our hypotheses do not imply that $Y$ is uniformly
simply-connected, so we must be a bit careful.
Before we delve into the detailed proof, we need two straightforward
lemmas based on the uniform $1$-con\-nec\-ted\-ness of $\E^2$.

\begin{lem}\label{classical one}
Let $\bar f:Y\to\E^2$ be a $(\lambda,C)$-\qy. There exists $\delta$
depending only on $\lambda$ and $C$ such that for any continuous map
$D_2:D^2\to Y$, there is a continous map $D_1:D^2\to \E^2$ such that
$d(\bar f(D_2(u)), D_1(u)) \le \delta$ for all $u\in D_2$.
\end{lem}

\begin{proof}
Choose a triangulation of $D^2$ such that the image by $D_2$ of each
$2$-simplex lies in a ball of radius $1$. The map $D_1$ is constructed
by induction over the skeleta of this triangulation.
\end{proof}

\begin{lem}[\protect\mbox{cf.~\cite[Lemma 8.6]{sm:qisurf}}]
\label{classical two}\mbox{}\\
For every $D>0$ there exists $\epsilon(D)$ such that any continuous
map $h:S^1\times \{0,1\} \to \E^2$ satisfying $d(h(t,0),h(t,1))\le D$
for all $t\in S^1$ can be extended to a continuous map $h: S^1\times
I\to \E^2$ such that $\diam(h(t\times I)) \le \epsilon(D)$ for all
$t\in S^1$.\cqfd
\end{lem}

\begin{proof}[Proof of Proposition~\ref{topo rigidity}]
Let $f:\E^2 \to Y$ and $\bar f:Y\to \E^2$ be
$(\lambda,C)$-quasi-isometries such that $d(\bar f(f(x)),x)\le C$ for
all $x\in \E^2$ and $d(f(\bar f(y)),y)\le C$ for all $y\in Y$.

Choose a point $x_0\in \E^2$ and let $c_1:S^1\to \E^2$ be an embedding
whose image is the circle around $x_0$ of radius $R$, where $R$ is a
large constant to be determined. The map $f\circ c_1$ avoids the ball
of radius $\lambda^{-1} R-C$ around $f(x_0)$. There is a continuous map
$c_2:S^1\to Y$ such that $d(f(c_1(t)),c_2(t)) \le 2(C+1)$ for all
$t$. Thus the image of $c_2$ avoids the ball of radius $\lambda^{-1}
R-3C-2$ around $f(x_0)$.

We want to fill $c_2$ with a continuous map $D_2:D^2\to Y$ which is
``far off''. Since $Y$ is \homeoic\ to $\Rr^3$, it is simply-connected
at infinity, so for any $R'\ge 0$ we can choose $R$ so that every loop
in $Y-B(f(x_0),\lambda^{-1}R-3C-2)$ can be filled in $Y-B(f(x_0),R')$.
We will see later how to choose $R'$, and therefore $R$.

Applying Lemma~\ref{classical one}, we get a constant
$\delta=\delta(\lambda,C)$ such that there is a continuous map
$D_1:D^2\to \E^2$ satisfying $d(D_1(u),\bar f(D_2(u))) \le \delta$ for
all $u$. We want to apply Lemma~\ref{classical two} with
$h(\cdot,0)=c_1$ and $h(\cdot,1)=\bord D_1$. Chasing through the
inequalities, we find that the hypothesis of this lemma is fulfilled
with $D=2(C+1)\lambda+2C+\delta$. Choose $R'$ large enough so that
$\epsilon(D) \le \lambda^{-1} R'-C-1$ and $\lambda^{-1} R'-C \ge
\delta +10$. Then Lemma~\ref{classical two} implies that $c_1$ and
$\bord D_1$ are homotopic in the complement of $x_0$. Furthermore,
$D_1$ misses $x_0$. Hence $c_1$ is null-homotopic in the complement of
$x_0$, which is a contradiction.
\end{proof}

\section{The 3-hyperbolic case}\label{sec cas hyperbole}
The goal of this section is to prove Theorem~\ref{prince} in the case
where $X_0$ is $3$-hyperbolic. In the first two subsections, we
develop some preliminary material. The proof itself is given
in subsection~\ref{proofitself}.

\subsection{Half-minima and Bloch principle}

\begin{defi}
  Let $X$ be a metric space and $h:X\to (0,+\infty)$ be a function. A
  \bydef{half-minimum} for $h$ is a point $x\in X$ such that
  $h(y)\ge\frac12h(x)$ for every $y\in B(x,\frac12\sqrt{h(x)})$.
\end{defi}

\begin{lem}[\protect{\cite[Lemma 7.3]{sm:qisurf}}]\label{demi minimum}
Let $X$ be a complete metric space and $h:X\to (0,+\infty)$ a function
which is locally bounded away from zero.
Let $x$ be a point of $X$ such that $h(x)<\frac12$. Then there exists
$x'\in B(x,2)$ such that $h(x') \le h(x)$ and $x'$ is a half-minimum. \cqfd
\end{lem}

The following result generalizes~\cite[Lemma~7.4]{sm:qisurf}.

\begin{theo}\label{bloch}
  Let $k\ge 2$ be an integer.  Let $(X,g)$ and $(X_0,g_0)$ be complete
  conformally flat \riem\ $k$-\var s of \bg.  Suppose that $(X_0,g_0)$
  has a cocompact group of isometries and is $k$-hyperbolic.
  Let $c:X\to X_0$ be a
  conformal embedding.  Define a function $\mu:X\to (0,+\infty)$ by setting
  $g=\mu^2 c^*g_\hyp$.  Then there is a constant $\mu_0>0$ such that
  $\mu(x)\ge\mu_0$ for all $x$.
\end{theo}

\begin{proof}
Seeking a contradiction, suppose that there is a sequence $x_n\in X$
such that $\mu(x_n)$ goes to $0$. By Lemma~\ref{demi minimum}
applied with $h=\mu$, there is no loss of generality in assuming that
each $x_n$ is a half-minimum.

Since $X$ is conformally flat and of \bg, there exist constants $r,\lambda$
and for each $n$ a conformal chart $\phi_n:B_{\E^k}(0,r)\to X$ such that
$\phi_n(0)=x_n$, $\|D\phi_n(0)\|=1$ and $\sup_{a\in B_{\E^k}(0,r)}
\|D\phi_n(a)\| \le 1/2\lambda$.

Set $B_n:=B_{\E^k}(0,\lambda/\sqrt{\mu(x_n)})$. Define
a mapping $z_n:B_n\to \E^k$ by $z_n(a):=\mu(x_n) a$. For large
$n$ we have $\lambda\sqrt{\mu(x_n)}\le r$, so the image of
$z_n$ lies in $B_{\E^k}(0,r)$. By hypothesis, we can for each $n$
postcompose $c$ with an isometry of $X_0$ so that the resulting map
$c_n:X\to X_0$ sends $x_n$ into a compact set $K$ independent of $n$. This
map $c_n$ is conformal.

Finally we set $f_n= c_n\circ \phi_n \circ z_n$. The goal is to find a
converging subsequence of $f_n$ and look at the limit to get a contradiction.

For this we need to estimate $\sup \|Df_n\|$ from above. First we see
that $\|Dz_n(a)\|=\mu(x_n)$ and $\|D\phi_n(z_n(a))\|\le
1/2\lambda$ for all $a\in B_n$. Thus, if $a\in B_n$, then $\phi_n(z_n(a))
\in B(x_n,\frac12\sqrt{\mu(x_n)})$ and the half-minimum property says
that $\mu(\phi_n(z_n(a)))\ge \frac12\mu(x_n)$. Now $c_n$ is conformal
with dilatation $1/\mu$, so $\|Dc_n(\phi_n(z_n(a)))\|\le
2/\mu(x_n)$. We deduce:

\begin{align*}
\|Df_n(a)\| &\le \|Dc_n(\phi_n(z_n(a)))\| \cdot
\|D\phi_n(z_n(a))\| \cdot \|Dz_n(a)\| \\
  &\le \frac2{\mu(x_n)} \cdot \frac1{2\lambda} \cdot \mu(x_n) \\
  &\le \frac1\lambda.
\end{align*}

As a consequence, for fixed $n$, the sequence $\{\tilde f_p\}_{p\ge
n}$ obtained by restricting each $f_p$ to $B_n$ is equicontinuous on
$B_n$, and for every $p$ we get $\tilde f_p(0)\in K$ and $\|D\tilde
f_n(0)\|=1$. By Ascoli's Theorem, $\{\tilde f_p\}_{p\ge n}$
subconverges. By diagonal extraction we get a sequence of conformal
mappings $g_n:B_n \to X_0$ which converges uniformly on compact
subsets to a mapping $g:\E^k \to X_0$. By general properties of
quasiconformal mappings, $g$ is $1$-quasiconformal or constant.

Now $X_0$ is $k$-hyperbolic and $\E^k$ is $k$-parabolic. Hence
by~\cite[Proposition 5.1]{chs:picard}, there is no quasiconformal
mapping from $\E^k$ to $X_0$. This implies that our map $g$ is constant.

Let us write $h_n$ for the restriction of $g_n$ to the unit ball. Since the
convergence is uniform on compact subsets, for large $n$ the image of
$h_n$ is contained in the image of a conformal chart $\phi$
for $X_0$. Now $\phi^{-1}\circ h_n$ is a sequence of conformal maps
between domains in $\E^k$. If $k=2$, such maps are holomorphic; if
$k\ge 3$ they are restrictions of Moebius transformations by
Liouville's Theorem. In any case
the condition $\|D\tilde f_n(0)\|=1$ gives a uniform lower bound on
the derivatives of $h_n$ at $0$, which leads to a contradiction.
\end{proof}

\subsection{Area and diameter estimates}\label{nouveau}

\begin{lem}\label{isodiam}
  Let $X$ be a \geod\ space \qi\ to a finitely \gen\ group. There
  exists a function $f_1:[0,+\infty) \to [0,+\infty)$ such that for
  every bounded subset $\Omega\subset X$ we have $\diam \Omega \le
  f_1(\diam \delta\Omega)$.
\end{lem}

\begin{proof}
In the Cayley graph of a finitely \gen\ group, each point is
contained in a (biinfinite) geodesic. Since $X$ is \qi\ to a finitely
\gen\ group, it has the corresponding ``quasi'' property: there exist
constants $\lambda,C\ge0$ such that for every $x\in X$, there is a
$(\lambda,C)$-\qg\ $\alpha:\Rr\to X$ such that $d(x,\alpha(0)) \le C$.
Since $X$ is \geod, we can assume without loss of generality that
$\alpha$ is continuous.

Let $\Omega\subset X$ be a bounded subset and $x$ be a point of
$\Omega$.  Consider a \qg\ $\alpha$ satisfying the above properties.
If $\alpha(0)\not\in \Int\Omega$, connect $x$ to $\alpha(0)$ by a
geodesic segment. This segment has to cross the frontier of $\Omega$,
so $d(x,\delta \Omega)\le C$.

Otherwise $\alpha$ meets $\delta\Omega$ for at least one negative time
$t_1$ and one positive time $t_2$. Without loss of generality, suppose
that $|t_1|\ge t_2$. Then $d(\alpha(0),\delta\Omega) \le \lambda t_2
+C$. Moreover,

\begin{align*}
t_2\le t_2-t_1&\le \lambda d(\alpha(t_1),\alpha(t_2)) +C\\
              &\le \lambda \diam(\delta\Omega) +C,
\end{align*}
hence
\begin{align*}
d(x,\delta\Omega)&\le C+d(\alpha(0),\delta\Omega)\\
              &\le \lambda (\lambda\diam(\delta\Omega) +C) +2C.
\end{align*}

Therefore, this inequality holds for each $x\in\Omega$. We conclude by setting
$f_1(r):=3(\lambda^2r+(\lambda+2)C)$.
\end{proof}

\begin{lem}\label{isodiam singulier}
  Let $X$ be a \riem\ $3$-\var\ \qi\ to a finitely generated group.
  Suppose that $H_2(X)=0$. Then there exists a function
  $f_2:[0,+\infty) \to [0,+\infty)$
  such that for every compact $K\subset X$ and every continuous map
  $s:S^2\to X-K$ that does not represent the trivial homology class in
  $X-K$, we have $\diam K \le f_2(\diam(s(S^2)))$.
\end{lem}

\begin{proof}
  To start off, thicken the image of $s$ into a codimension 0 \svar\ $Y$
  contained in $X-K$ and such that $\diam Y\le \diam (s(S^2))+1$.
  Then represent the class of $s$ in $H_2(Y)$ by a system of embedded
  surfaces. Since $[s]\neq0 \in H_2(X-K)$, one of these surfaces, say
  $F$, is not trivial in $H_2(X-K)$. Since $H_2(X)=0$, $F$ bounds a
  compact \svar\ $\Omega$.

Now $K\subset\Omega$ and $\diam(\delta\Omega)=\diam(F)\le \diam
s(S^2)+1$, so we conclude by applying Lemma~\ref{isodiam}.
\end{proof}

\begin{prop}\label{coupe sphere}
For all $A,\epsilon>0$ there exists $L=L(A,\epsilon)$ such that
if $S$ is a \riem\ $2$-sphere of area at most $A$
and $\gamma\subset S$ is an embedded curve, then there is a system
$\{\xi_1,\ldots,\xi_n\}$ of pairwise disjoint embedded curves that cobound
a \svar\ $U$ such that $\gamma\subset U\subset N(\gamma,L)$ and
$|\xi_i|\le \epsilon$ for all $i$.
\end{prop}

\begin{proof}
Without loss of generality assume that $|\gamma|\ge\epsilon$.  By
Loewner's Theorem (see e.g.~\cite{glp:struc}), there is a constant
$L'$ such that any \riem\ annulus of area at most $2A$ and whose
boundary components are a distance at least $L'$ has systole at most
$\epsilon$ (recall that the \bydef{systole} is the infimum of
lengths of noncontractible curves). Set $L:=L'+\epsilon$.
  
Let us give some definitions and notation.  Let $\Xi$ be a system of curves
embedded in $N(\gamma,L)-\gamma$. We denote by $E(\Xi)$ the set
of points $x\in S-N(\gamma,L)$ such that there exists $\xi\in \Xi$
which separates $x$ from $\gamma$. We shall say that $\Xi$ is
\bydef{embedded} if its elements are pairwise disjoint.

The conclusion of Propposition~\ref{coupe sphere} can be reformulated
as follows: there exists an embedded system $\Xi$ satisfying
conditions (i) and (ii) below:

\begin{enumerate}
\item $|\xi| \le \epsilon$ for all $\xi\in\Xi$~;
\item  $E(\Xi) = S-N(\gamma,L)$.
\end{enumerate}

Here is our strategy: first we find a (possibly not embedded) system
satisfying (i) and (ii). Then if it is not embedded, we show how to do
surgery on it to produce a system with fewer self-intersections still
satisfying (i) and (ii).

\begin{lem}\label{pas plonge}
  There is a system $\Xi_0=\{\xi'_0,\ldots,\xi'_m\}$ satisfying
  conditions (i) and (ii).
\end{lem}

\begin{proof}
  For a generic choice of $L'$, $N(\gamma,L')$ is a planar surface.
  Let $\eta_0,\ldots,\eta_m$ be its boundary components. Let
  $D_1,\ldots, D_0$ be discs such that $\bord D_i=\eta_i$ and $D_i\cap
  \gamma= \emptyset$.
  
  Fix $i$ between $0$ and $m$. Let $Y_i$ be the annulus cobounded by
  $\gamma$ and $\eta_i$. Since $d(\gamma,\eta_i)=L'$, there is an
  essential curve $\xi'_i\subset Y_i$ of length at most $\epsilon$.
  Since $\xi'_i$ is essential, it cannot lie in one of the $D_j$'s.
  Hence it is contained in $N(\gamma,L)$. This ensures that the system
  $\{\xi'_i\}$ satisfies (ii). By construction it also satisfies (i).
\end{proof}

Before proceeding, we need one more piece of notation. Let $\Xi$ be an
embedded system of curves and $\eta$ an embedded curve in $S-\gamma$
in general position with respect to $\Xi$. Let $\sing(\Xi,\eta)$
denote the cardinal of $\eta \cap \Xi$. In particular,
$\sing(\Xi,\eta)=0$ iff $\Xi \cup \{\eta\}$ is embedded.

Assume that $\sing(\Xi,\eta)>0$. A well-known lemma (cf.~\cite{flp})
ensures that there is a \bydef{bigon} between $\eta$ and some
$\xi\in\Xi$, i.e.~a disc $D\subset S-\gamma$ whose boundary is the
union of two arcs $\alpha_1,\alpha_2$ with $\alpha_1\subset\eta$,
$\alpha_2\subset \xi$, and such that $\Int D$ meets neither $\eta$ nor
$\Xi$. We say that $\alpha_1$ (resp.~$\alpha_2$) is \bydef{exterior} to $D$
if the curve $\eta$ (resp.~$\xi$) surrounds $D$ (i.e.~$D$ is contained
in the unique disc bounded by this curve). Otherwise we say $\alpha_1$
(resp.~$\alpha_2$) is \bydef{interior}.

\begin{figure}[ht]
\begin{center}
\input{bigons.pstex_t}
\end{center}
\caption {\label{fig:bigons} Various types of bigons.}
\end{figure}
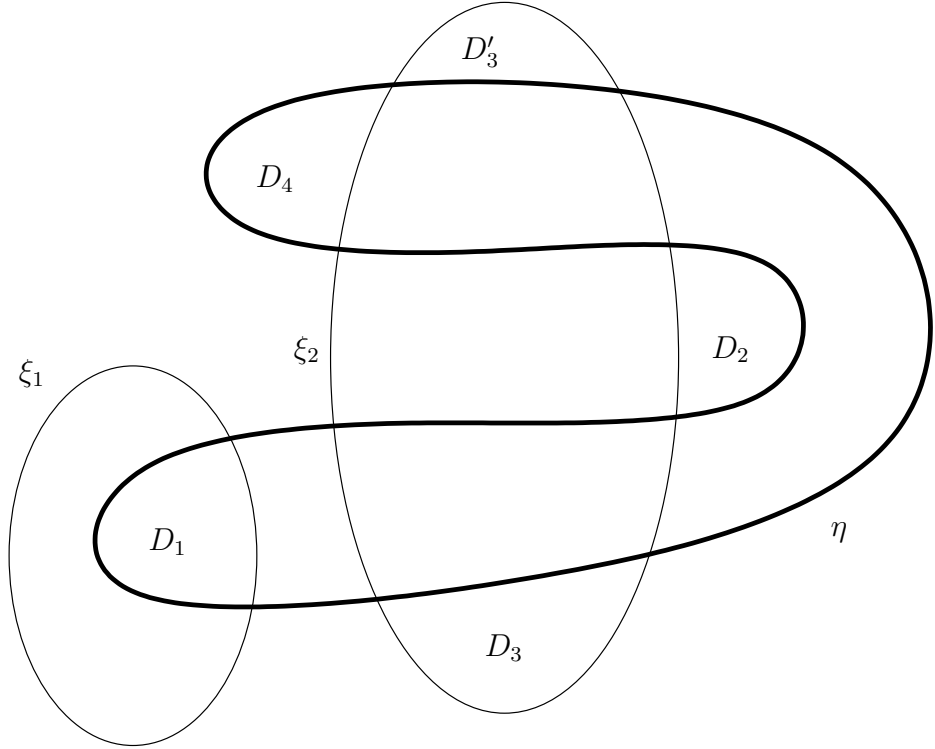

We will always use this notation, i.e.~if the bigon is called $D$, we
call $\alpha_1$ the arc lying in $\eta$ and $\alpha_2$ the arc lying
in a element of $\Xi$.  This allows to distinguish several types of
bigons: a bigon $D$ is called
\bydef{exterior} if (with the same notation as above) both arcs
$\alpha_1,\alpha_2$ are interior to $D$, \bydef{interior} if
$\alpha_1,\alpha_2$ are exterior to $D$, and \bydef{mixed} otherwise.
Furthermore, a mixed bigon $D$ is of \bydef{type~1} if $\alpha_1$ is
the interior arc, and of \bydef{type~2} otherwise.  Two mixed bigons
are \bydef{paired} if they are both of type~1 and involve the same
element of $\Xi$.

In the example illustrated by Figure~\ref{fig:bigons}, the curve
$\eta$ has five bigons of intersection with $\Xi=\{\xi_1,\xi_2\}$;
$D_1$ is interior, $D_2$ is exterior, $D_3,D'_3$ form a pair of mixed
bigons of type~1, and $D_4$ is mixed of type~2.

\begin{lem}\label{desingularise}
Let $\Xi$ be an embedded system and $\eta$ be an embedded curve in
$S-\gamma$ in general position with respect to $\Xi$. Suppose that
$\Xi \cup \{\eta\}$ satisfies condition~(i). Then there is an embedded
system $\Xi'$ satisfying~(i) and such that $E(\Xi') \supset
E(\Xi) \cup E(\eta)$.
\end{lem}

\begin{proof}
The proof is by induction on $\sing(\Xi,\eta)$. If this number is
$0$, we can just set $\Xi':=\Xi \cup \{\eta\}$. Otherwise we will
show how to use the induction hypothesis by applying to $(\Xi,
\eta)$ one or more of the three operations described below.

\paragraph{($\mathrm{T}_0$)}
Let $\xi$ be an element of $\Xi$ lying in the interior of a disc
$D\subset S-\gamma$ bounded by $\eta$ or by another element of $\Xi$.
The operation $\mathrm{T}_0$ consists in keeping $\eta$ unchanged and
removing $\xi$ of $\Xi$.

We say that $(\Xi,\eta)$ is \bydef{reduced} if $\mathrm{T}_0$ cannot
be applied to it. Clearly, $\mathrm{T}_0$ does not change
condition~(i) nor $E(\Xi) \cup E(\eta)$ and never increases
$\sing(\Xi,\eta)$. Hence we can always assume that $(\Xi,\eta)$ is
reduced.

\paragraph{($\mathrm{T}_1$)}
Let $D$ be a bigon bounded by arcs $\alpha_1\subset\eta$ and
$\alpha_2\subset \xi$. Let $\beta$ be the closure of $\eta-\alpha_1$.
Replace $\eta$ by a curve obtained from $\beta\cup \alpha_2$ by a
small isotopy that removes the intersections in the neighborhood of 
$\alpha_2$. We call this operation \bydef{pushing $\eta$ through $D$}.

Operation $\mathrm{T}_1$ decreases $\sing(\Xi,\eta)$ by $2$. If
$|\alpha_1|\ge |\alpha_2|$, it respects condition~(i) (taking the
isotopy sufficiently small); furthermore, $E(\Xi) \cup E(\eta)$ can go
down only if the bigon is mixed and $\alpha_1$ is the exterior arc.

Symmetrically we define \bydef{pushing $\xi$ through $D$}.

\paragraph{($\mathrm{T}_2$)}
Add to $\Xi$ a curve disjoint from $\Xi$, contained in a bigon $D$ and
obtained from $\bord D$ by a small isotopy.

Note that $E(\Xi) \cup E(\eta)$ never goes down when $\mathrm{T}_2$ is
applied. Condition~(i) is preserved provided that $|\bord D| \le \epsilon$.

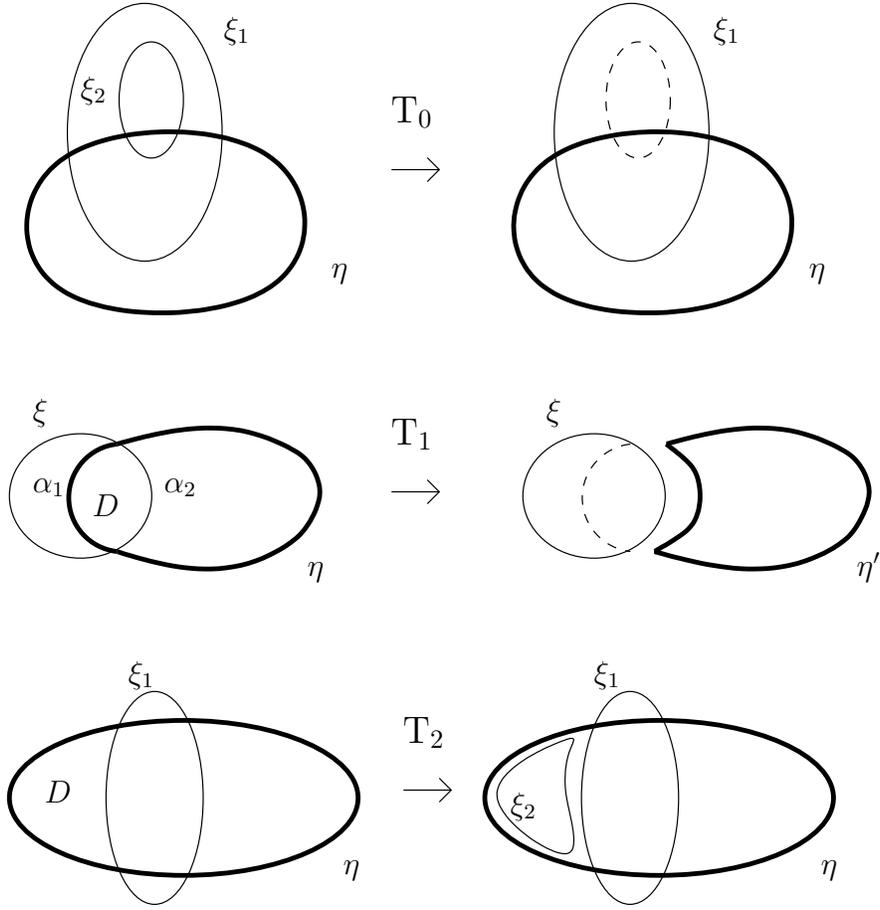
\begin{figure}[hbt]
\begin{center}
\input{moves.pstex_t}
\end{center}
\caption{\label{fig:moves} The three moves $\mathrm{T}_0$--$\mathrm{T}_2$.}
\end{figure}

To deal with mixed bigons, we need the following lemma:
\begin{slem}\label{oppose}
  Let $\Xi$ be an embedded system of curves and $\eta$ an embedded
  curve in general position with respect to $\Xi$. Suppose that
  $(\Xi,\eta)$ is reduced and there are no interior bigons. Then
  either $(\Xi,\eta)$ is embedded, or there are paired bigons.
\end{slem}

\begin{proof}
Assume that $(\Xi,\eta)$ is not embedded. Choose $\xi\in \Xi$ such
that $\xi\cap\eta\neq\emptyset$. Let $D_\xi$ be the disc bounded by
$\xi$. Then $\eta\cap D_\xi$ consists of one or more arcs. Let $\beta$ be
one of these arcs. On each side of $\beta$ choose an outermost arc
$\beta_1$ (resp.~$\beta_2$). Then each $\beta_i$ cobounds a bigon
$D_i$ with a subarc of $\xi$. Each $D_i$ is interior to $\xi$. Since
by hypothesis there are no interior bigons, $D_1$ and $D_2$ must be
mixed of type 1, hence paired.
\end{proof}

We turn to the proof of Lemma~\ref{desingularise}. Take a non-embedded
pair $(\Xi,\eta)$ fulfilling the hypotheses of this lemma and assume
the result holds for all pairs $(\Xi',\eta')$ such that
$\sing(\Xi',\eta') \le \sing(\Xi,\eta)$. By applying move
$\mathrm{T}_0$ zero or more times, we can assume that $(\Xi,\eta)$ is
reduced. If there is an exterior bigon or an interior bigon, we can
perform move $\mathrm{T}_1$, pushing the bigger arc through the smaller
one, and apply the induction hypothesis to the resulting configuration.

If all bigons are mixed, then by Sublemma~\ref{oppose} there are some
paired bigons. Let us write $\alpha_1 \cup \alpha_2$ and
$\alpha'_1 \cup \alpha'_2$ for their boundaries, with the usual conventions.
We have $|\alpha_1|+|\alpha'_1|\le \epsilon$
and $|\alpha_2|+|\alpha'_2|\le \epsilon$. Hence
$|\alpha_1|+|\alpha_2|+|\alpha'_1|+|\alpha'_2|\le 2\epsilon$, and we
may assume without loss of generality that $|\alpha_1|+|\alpha_2|\le
\epsilon$ and $|\alpha_1| \ge |\alpha_2|$. Then we apply
$\mathrm{T}_2$ followed by $\mathrm{T}_1$. This finishes the proof of
Lemma~\ref{desingularise}.
\end{proof}

At last we prove Proposition~\ref{coupe sphere}. Consider the
system $\Xi_0=\{\xi'_0,\ldots,\xi'_m\}$ given by Lemma~\ref{pas
plonge}. Applying Lemma~\ref{desingularise} with $\Xi=\{\xi'_0\}$ and
$\eta=\xi'_1$, we get an embedded system $\Xi_1$ satisfying~(i) and
such that $E(\Xi_1) \supset E(\{\xi'_0,\xi'_1\})$. Then we apply
Lemma~\ref{desingularise} successively for each $i$ from $2$ to $m$,
putting $\Xi=\Xi_{i-1}$ and $\eta=\xi'_i$. The outcome of each step is
an embedded system $\Xi_i$ satisfying~(i) such that $E(\Xi_1) \supset
E(\{\xi'_0,\ldots,\xi'_i\})$. Hence $\Xi_m$ is embedded and
satisfies~(i) and (ii). This completes the proof of Lemma~\ref{coupe sphere}.
\end{proof}

\subsection{Proof of the 3-hyperbolic case}\label{proofitself}
Let $(X,g)$ and $(X_0,g_0)$ be \riem\ $3$-\var s satisfying the
hypotheses of Theorem~\ref{prince}. Assume that $X_0$ (and hence $X$)
is $3$-hyperbolic. Let $c:(X,g)\to (X_0,g_0)$ be a conformal \diffeo.
Let $\mu:X\to (0,+\infty)$ denote the function defined by $g=\mu^2
c^*g_0$. We sometimes still denote (abusively) by $\mu$ the function
$\mu \circ c^{-1}$. We shall prove:

\begin{prop}\label{autre sens}
There are constants $r_0,\mu_1>0$ such that for every $x\in X_0$ the
following holds:
$$\diam c^{-1}(B(x,r_0/2))\le \mu_1.$$
\end{prop}

To see why this implies Theorem~\ref{prince}, we note that
Theorem~\ref{bloch} implies that $d(x,y)\ge \mu_0 d(c(x),c(y))$ for
some constant $\mu_0>0$ and all $x,y\in X$. Since our goal is to prove
that $X$ is \qi\ to $X_0$, we need an upper bound for $d(x,y)$ in
terms of $d(c(x),c(y))$. For all $x,y\in X$, choose a geodesic arc
$\gamma$ connecting $c(x)$ to $c(y)$. We can cover $\gamma$ by $n$ balls
of radius $r_0/2$ with $n$ bounded above by a linear function of
$d(c(x),c(y))$. Hence by Proposition~\ref{autre sens}, $d(x,y)$ is bounded
by a linear function of $d(c(x),c(y))$. This proves that $c$ is a \qy,
which completes the proof of Theorem~\ref{prince} in the
$3$-hyperbolic case.\cqfd

The remainder of this paper is devoted to the proof of
Proposition~\ref{autre sens}. First of all, we gather in the next
lemma some immediate consequences of the bounded
geometry hypothesis on $X_0$. 
 \begin{lem}\label{consequences bg}
There exist positive constants $R_0,\lambda_1,\lambda_2,\lambda_3,N_1$
such that for every $x\in X_0$ and every $\rho \le R_0$ we have:
\begin{enumerate}
\item $|S(x,\rho)|\le \lambda_1^2\rho^2$.
\item For all $x_1,\ldots,x_n\in S(x,\rho)$ such that if $i\neq j$, then
$d(x_i,x_j) \ge 9\rho/5$, we have $n\le N_1$.
\item For each parallel circle $\gamma$ on $S(x,\rho)$ of latitude at most
$\pi/5$, we have $\diam(\gamma)\ge\lambda_2\rho$.
\item Let $\xi$ be a curve on $S(x,\rho)$ of length bounded above by
$\lambda_3\rho$. Then exactly one of the two discs bounded by $\xi$ on
$S(x,\rho)$ has diameter bounded above by $3\lambda_2\rho/4$.
\end{enumerate}
\end{lem}

\begin{proof}
Choose for $R_0$ a lower bound for the injectivity radius of
$X_0$. The restriction of the exponential map at any point $x$ to the
ball of radius $R_0$ around the origin is a bilipschitz embedding with
uniform Lipschitz constant.

In assertion~(iii), the word ``parallel'' refers to the image by the
exponential map at $x$ of a parallel for the standard spherical
coordinates in $\Rr^3$; the word ``latitude'' is to be interpreted in
the same sense. Since $R_0$ is less than the injectivity radius at $x$,
$\gamma$ is indeed a topological circle.
\end{proof}

\begin{lem}\label{majore L}
\begin{enumerate}
\item There is a constant $\lambda_4>0$ such that for every $\nu>0$,
every $x\in X_0$ and every $r\le R_0$, if $|c^{-1}(B(x,r))|\le \nu$,
then there exists $\rho\in[9r/{10},r]$ such that
$|c^{-1}(S(x,\rho))|\le \lambda_4 \nu^{2/3}$.
\item For every $A>0 $ there exists $L_1=L_1(A)$ such that for every
$x\in X_0$ and every $\rho\le R_0$, if $|c^{-1}S(x,\rho)| \le A$, then there
is a curve $\gamma\subset c^{-1}S(x,\rho)$ satisfying $|\gamma|
\le L_1$ and $\diam c(\gamma)\ge \lambda_2\rho$.
\end{enumerate}
\end{lem}

\begin{proof}
Set $$\lambda_4 := {\left(\frac{\lambda_1} {\ln(10/9)}\right)}^{2/3} +1.$$

If (i) does not hold, then by H\"older's inequality the following is
true for all $\rho\in [9r/{10},r]$:
\begin{align*}
\lambda_4^{3/2} \nu &\le \left(\int_{S(x,\rho)}\mu^2\,d\vol\right)^{3/2}
\le |S(x,\rho)|^{1/2}\cdot \int_{S(x,\rho)}\mu^3\,d\vol\\
&\le \lambda_1\rho \cdot \int_{S(x,\rho)}\mu^3\,d\vol.
\end{align*}

Dividing by $\rho$ and integrating between $9r/{10}$ and $r$, we get:
$$\lambda_4^{3/2} \nu \cdot \int_{9r/{10}}^r \frac{d\rho}\rho \le \lambda_1
|c^{-1}B(x,r))|\le \lambda_1 \nu,$$
so $$\lambda_4^{3/2} \ln({10}/9)\le \lambda_1,$$
which contradicts the choice of $\lambda_4$.

The proof of (ii) is similar, using~\ref{consequences bg}(iii) and the
Cauchy-Schwarz inequality in spherical coordinates.
\end{proof}

Given $x\in X_0$ and $\nu>0$, we let $r(x,\nu)$ denote the infimum of
 numbers $\rho>0$ such that $|c^{-1} (B(x,\rho))|\ge \nu$. For fixed
 $\nu$, the function $x\mapsto r(x,\nu)$ may not be continuous, but it
 is locally bounded away from zero, so we can apply
 Lemma~\ref{demi minimum} to it.

\begin{lem}\label{nu zero}
There is a constant $\nu_0>0$ such that $\inf\{r(x,\nu_0)\mid
x\in X_0\}>0$.
\end{lem}

\begin{proof}
First we reduce this lemma to the following claim:

\begin{claim}
For every $\nu>0$, if
$\inf \{r(x,\nu)\mid x\in X_0\}=0$, then there is a domain
$\Omega\subset X$ such that $|\Omega| \ge \nu$ and $|\bord\Omega|
\le N_1\lambda_4 \nu^{2/3}$.
\end{claim}

Let us prove by contradiction that this claim implies Lemma~\ref{nu zero}.
Let $\nu_i\to +\infty$ be a sequence such that $\inf\{r(x,\nu_i)\mid
x\in X_0\}=0$. The claim supplies a sequence of domains $\Omega_i\subset X$
satisfying $|\Omega_i|\to +\infty$ and
$$\frac{|\bord\Omega_i|}{|\Omega_i|^{2/3}}\le N_1\lambda_4.$$

It follows that for all $D>3$,
$$\frac{|\bord\Omega_i|}{|\Omega_i|^{(D-1)/D}}\to 0,$$
which shows that $X$ has \isop\ dimension at most $3$. By
Theorem~\ref{varogromov},  $X$ is
$3$-parabolic. This contradicts the conformal invariance of $3$-parabolicity.

The next task is to prove the claim. Fix $\nu>0$ and apply
Lemma~\ref{demi minimum} to $x\mapsto r(x,\nu)$. This gives a point
$x_\nu\in X_0$ satisfying $r(x_\nu,\nu)<R_0$ and such that for any
$x\in X_0$, if $d(x,x_\nu)\le\frac12\sqrt{r(x_\nu,\nu)}$ then
$r(x,\nu)\ge\frac12 r(x_\nu,\nu)$. For simplicity, let us write $r_\nu$
for $r(x_\nu,\nu)$. Without loss of generality assume that
$r_\nu<\frac14\sqrt{r_\nu}$. Then for every $x\in B(x_\nu,2r_\nu)$ and
every $r\le r_\nu$ the inequality $|c^{-1}(B(x,r))|\le \nu$ holds.

By Lemma~\ref{majore L}(i) applied with $r=r_\nu$, we can associate to
each point $x\in X_0$ such that $d(x,x_\nu)=r_\nu$ a number
$\rho(x)\in[9r_\nu/{10},r_\nu]$ satisfying $|c^{-1}(\bord B(x,\rho(x)))|
\le \lambda_4 \nu^{2/3}$. Let $\{x_1,\ldots,x_n\}$ be a set of points
of $S(x_\nu,r_\nu)$ with $n$ minimal such that the metric balls
$B(x_i,\rho(x_i))$ cover $S(x_\nu,r_\nu)$. Set $\Omega
:= c^{-1}(B(x_\nu,r_\nu)\cup \bigcup_i B(x_i,\rho(x_i)))$. Then
$\Omega$ contains $c^{-1}(B(x_\nu,r_\nu))$, so $|\Omega|\ge
\nu$. Since $n$ is minimal, the balls $B(x_i,\rho(x_i))/3$ are
pairwise disjoint. Since $\rho(x_i)\ge 9r_\nu/{10}$ and $r_\nu<R_0$,
Lemma~\ref{consequences bg}(ii) gives $n\le N_1$. It follows that 
$|\bord\Omega| \le n \lambda_4 \nu^{2/3}\le N_1\lambda_4 \nu^{2/3}$.
\end{proof}

Set $r_0:=\min(R_0,\inf\{r(x,\nu_0) \mid x\in X_0\})$.
Applying~\ref{majore L}(i) with $\nu=\nu_0$ and $r=r_0$, we obtain
a function $\rho:X_0\to [9r_0/{10},r_0]$ satisfying
$$|c^{-1}(S(x,\rho(x)))|\le \lambda_4 \nu_0^{2/3}.$$

Define $A_0:=\lambda_4 \nu_0^{2/3}$.
Fix a point $x\in X_0$ and consider the metric sphere $S=S(x,\rho(x))$.
By construction, $c^{-1}(S)$ is not null-homotopic in the complement of
$c^{-1}(B(x,r_0/2))$; furthermore, its area is bounded above by $A_0$.
If we had a uniform upper bound of the diameter of $S$ (as opposed to
the area) we could apply Lemma~\ref{isodiam}. Since we do not have such a
bound, we are going to use Proposition~\ref{coupe sphere} and
Lemma~\ref{isodiam singulier}. 

Lemma~\ref{majore L}(ii) provides a constant $L_1=L_1(A_0)$ and an
embedded curve 
$\gamma\subset c^{-1}(S)$ such that $|\gamma|\le L_1$ and $\diam
c(\gamma)\ge \lambda_2\rho(x)$. Let us apply
Lemma~\ref{coupe sphere} with $A=A_0$ and $\epsilon=\inf\{\inj(X), 
9\mu_0\lambda_3 r_0/{10}\}$. Each $\xi_i$ bounds a small disc $D_i$ of
diameter at most $9\mu_0\lambda_3 r_0/{10}$. It follows that
$\diam c(D_i)\le 9\lambda_3 r_0/{10}\le \lambda_3\rho(x)$. Therefore, of
the two discs bounded by $c(\xi)$ on $S$, the small one
is homotopic (with fixed boundary) to $c(D_i)$ in the complement of
$B(x,r_0/2)$.

Hence $U\cup \bigcup_i D_i$ is a (possibly not embedded)
$2$-sphere homotopic to $c^{-1}(S)$, of diameter at most
$C_4+2(L+9\mu_0\pi r_0/{100})$. To conclude, we apply
Lemma~\ref{isodiam singulier} and set $\mu_1:=f_2(C_4+2(L+9\mu_0\pi
r_0/{100}))$. The proof of Proposition~\ref{autre sens}, and hence of
Theorem~\ref{prince}, is now complete.

\section{Final remarks}\label{sec final}

We already remarked that in the $3$-parabolic case we proved a
stronger statement (Theorem~\ref{parabolic case}). In particular we
did not use the hypothesis of conformal flatness. Likewise, in the
$3$-hyperbolic case we did not use the hypothesis that
$X_0\cong\Rr^3$. Hence we have actually proved:

\begin{theo}\label{strong hyperbolic}
Every complete, conformally flat, $3$-hyperbolic Riemannian $3$-\var\
which admits a geometric group action is large-scale conformally rigid.
\end{theo}

In another direction, the methods of this paper can be used to show:

\begin{theo}\label{dim two}
Every complete Riemannian $2$-\var\ which admits a geometric group
action is large-scale conformally rigid.
\end{theo}

For $\E^2$ and $\H^2$ this result is essentially proved in~\cite{mess:seifert}.
The method used there for the hyperbolic plane does not seem to extend to
nonsimply-connected hyperbolic surfaces. Partial results along the
same lines (using other techniques closer to those of the present paper)
were obtained in~\cite{sm:qisurf}.

One may ask whether our hypothesis that all \var s involved
be of \bg\ is necessary.
Our definition of \bg\ is stronger than those usually found in the
literature, i.e.~uniform bounds on Ricci curvature, sectional
curvature, and/or injectivity radius. Some results (e.g.~those of
section~\ref{sec parabole}) hold under weaker assumptions. 
Note however that any \var\ with sectional curvature
bounded in absolute value and injectivity radius bounded away from
zero is \qi\ to a \var\ of \bg\ in our sense. To see this, construct a
triangulation with controlled geometry, take the regular piecewise
Euclidean metric associated to this triangulation, and smoothen out
consistently the singularities (cf.~\cite{sm:qisurf,attie:bg}). For
this reason, issues of minimal hypotheses were ignored in this paper.

Theorem~\ref{strong hyperbolic} applies to all $3$-hyperbolic \var s
which are regular covers of closed $3$-\var s that admit conformally
flat structures. This is a large class of $3$-\var s, which includes
e.g.~lots of open hyperbolic \var s, but also many
others~(see~\cite{kapovich:survey} and the references therein).

Consider Thurston's eight geometries. By Theorem~\ref{varogromov}, the
$3$-parabolic ones are $\Ss^3$, $\Ss^2\times\Rr$, and $\E^3$. They are
all large-scale conformally rigid (this is trivial for the first two,
and follows from~Theorem~\ref{prince} for the third one). Our main
theorem applies to $\H^3$ and $\H^2\times \Rr$. It does not apply to
$\slr$, but since there exist closed $\slr$-\var s with conformally flat
metrics, any group \qi\ to $\slr$ is also \qi\ to some complete
conformally flat \var\ to which Theorem~\ref{prince} applies. This
raises an obvious question: is the large-scale conformal rigidity
property invariant by \qy?

Our result does not give anything for fundamental groups of closed
$\mathbf{Nil}$ and $\mathbf{Sol}$ manifolds, since they do not
admit flat conformal structures~\cite{goldman:noflat}. One can ask
whether $\mathbf{Nil}$ and $\mathbf{Sol}$ are large-scale conformally
rigid. The arguments of the present paper seem to suggest a positive answer.

\def\cprime{$'$}

\end{document}

%% file: bigons.pstex_t
\begin{picture}(0,0)%
\special{psfile=bigons.pstex}%
\end{picture}%
\setlength{\unitlength}{3947sp}%
\begingroup\makeatletter\ifx\SetFigFont\undefined%
\gdef\SetFigFont#1#2#3#4#5{%
  \reset@font\fontsize{#1}{#2pt}%
  \fontfamily{#3}\fontseries{#4}\fontshape{#5}%
  \selectfont}%
\fi\endgroup%
\begin{picture}(5828,4685)(3224,-5597)
\put(4114,-4345){\makebox(0,0)[lb]{\smash{\SetFigFont{12}{14.4}{\familydefault}{\mddefault}{\updefault}{\color[rgb]{0,0,0}$D_1$}%
}}}
\put(4788,-2062){\makebox(0,0)[lb]{\smash{\SetFigFont{12}{14.4}{\familydefault}{\mddefault}{\updefault}{\color[rgb]{0,0,0}$D_4$}%
}}}
\put(3301,-3286){\makebox(0,0)[lb]{\smash{\SetFigFont{12}{14.4}{\familydefault}{\mddefault}{\updefault}{\color[rgb]{0,0,0}$\xi_1$}%
}}}
\put(5026,-3136){\makebox(0,0)[lb]{\smash{\SetFigFont{12}{14.4}{\familydefault}{\mddefault}{\updefault}{\color[rgb]{0,0,0}$\xi_2$}%
}}}
\put(8401,-4261){\makebox(0,0)[lb]{\smash{\SetFigFont{12}{14.4}{\familydefault}{\mddefault}{\updefault}{\color[rgb]{0,0,0}$\eta$}%
}}}
\put(6226,-5011){\makebox(0,0)[lb]{\smash{\SetFigFont{12}{14.4}{\familydefault}{\mddefault}{\updefault}{\color[rgb]{0,0,0}$D_3$}%
}}}
\put(7651,-3136){\makebox(0,0)[lb]{\smash{\SetFigFont{12}{14.4}{\familydefault}{\mddefault}{\updefault}{\color[rgb]{0,0,0}$D_2$}%
}}}
\put(6076,-1261){\makebox(0,0)[lb]{\smash{\SetFigFont{12}{14.4}{\familydefault}{\mddefault}{\updefault}{\color[rgb]{0,0,0}$D'_3$}%
}}}
\end{picture}

%% file: moves.pstex_t
\begin{picture}(0,0)%
\special{psfile=moves.pstex}%
\end{picture}%
\setlength{\unitlength}{3947sp}%
\begingroup\makeatletter\ifx\SetFigFont\undefined%
\gdef\SetFigFont#1#2#3#4#5{%
  \reset@font\fontsize{#1}{#2pt}%
  \fontfamily{#3}\fontseries{#4}\fontshape{#5}%
  \selectfont}%
\fi\endgroup%
\begin{picture}(5456,5676)(3503,-5961)
\put(5926,-1036){\makebox(0,0)[lb]{\smash{\SetFigFont{14}{16.8}{\familydefault}{\mddefault}{\updefault}{\color[rgb]{0,0,0}$\mathrm{T}_0$}%
}}}
\put(5926,-3061){\makebox(0,0)[lb]{\smash{\SetFigFont{14}{16.8}{\familydefault}{\mddefault}{\updefault}{\color[rgb]{0,0,0}$\mathrm{T}_1$}%
}}}
\put(6001,-4936){\makebox(0,0)[lb]{\smash{\SetFigFont{14}{16.8}{\familydefault}{\mddefault}{\updefault}{\color[rgb]{0,0,0}$\mathrm{T}_2$}%
}}}
\put(8551,-2011){\makebox(0,0)[lb]{\smash{\SetFigFont{12}{14.4}{\familydefault}{\mddefault}{\updefault}{\color[rgb]{0,0,0}$\eta$}%
}}}
\put(5626,-5761){\makebox(0,0)[lb]{\smash{\SetFigFont{12}{14.4}{\familydefault}{\mddefault}{\updefault}{\color[rgb]{0,0,0}$\eta$}%
}}}
\put(8626,-5761){\makebox(0,0)[lb]{\smash{\SetFigFont{12}{14.4}{\familydefault}{\mddefault}{\updefault}{\color[rgb]{0,0,0}$\eta$}%
}}}
\put(4276,-4561){\makebox(0,0)[lb]{\smash{\SetFigFont{12}{14.4}{\familydefault}{\mddefault}{\updefault}{\color[rgb]{0,0,0}$\xi_1$}%
}}}
\put(7201,-4561){\makebox(0,0)[lb]{\smash{\SetFigFont{12}{14.4}{\familydefault}{\mddefault}{\updefault}{\color[rgb]{0,0,0}$\xi_1$}%
}}}
\put(6901,-2911){\makebox(0,0)[lb]{\smash{\SetFigFont{12}{14.4}{\familydefault}{\mddefault}{\updefault}{\color[rgb]{0,0,0}$\xi$}%
}}}
\put(3676,-2911){\makebox(0,0)[lb]{\smash{\SetFigFont{12}{14.4}{\familydefault}{\mddefault}{\updefault}{\color[rgb]{0,0,0}$\xi$}%
}}}
\put(4051,-3511){\makebox(0,0)[lb]{\smash{\SetFigFont{12}{14.4}{\familydefault}{\mddefault}{\updefault}{\color[rgb]{0,0,0}$D$}%
}}}
\put(3751,-5311){\makebox(0,0)[lb]{\smash{\SetFigFont{12}{14.4}{\familydefault}{\mddefault}{\updefault}{\color[rgb]{0,0,0}$D$}%
}}}
\put(6676,-5386){\makebox(0,0)[lb]{\smash{\SetFigFont{12}{14.4}{\familydefault}{\mddefault}{\updefault}{\color[rgb]{0,0,0}$\xi_2$}%
}}}
\put(3676,-3361){\makebox(0,0)[lb]{\smash{\SetFigFont{12}{14.4}{\familydefault}{\mddefault}{\updefault}{\color[rgb]{0,0,0}$\alpha_1$}%
}}}
\put(4876,-511){\makebox(0,0)[lb]{\smash{\SetFigFont{12}{14.4}{\familydefault}{\mddefault}{\updefault}{\color[rgb]{0,0,0}$\xi_1$}%
}}}
\put(7951,-511){\makebox(0,0)[lb]{\smash{\SetFigFont{12}{14.4}{\familydefault}{\mddefault}{\updefault}{\color[rgb]{0,0,0}$\xi_1$}%
}}}
\put(3976,-886){\makebox(0,0)[lb]{\smash{\SetFigFont{12}{14.4}{\familydefault}{\mddefault}{\updefault}{\color[rgb]{0,0,0}$\xi_2$}%
}}}
\put(5401,-3886){\makebox(0,0)[lb]{\smash{\SetFigFont{12}{14.4}{\familydefault}{\mddefault}{\updefault}{\color[rgb]{0,0,0}$\eta$}%
}}}
\put(5551,-2011){\makebox(0,0)[lb]{\smash{\SetFigFont{12}{14.4}{\familydefault}{\mddefault}{\updefault}{\color[rgb]{0,0,0}$\eta$}%
}}}
\put(4501,-3361){\makebox(0,0)[lb]{\smash{\SetFigFont{12}{14.4}{\familydefault}{\mddefault}{\updefault}{\color[rgb]{0,0,0}$\alpha_2$}%
}}}
\put(8851,-3886){\makebox(0,0)[lb]{\smash{\SetFigFont{12}{14.4}{\familydefault}{\mddefault}{\updefault}{\color[rgb]{0,0,0}$\eta'$}%
}}}
\end{picture}